\documentclass[leqno,12pt]{article}
\usepackage{amsmath}
\usepackage{amsfonts}
\usepackage{amssymb}
\usepackage{graphicx}
\newtheorem{theorem}{Theorem}[section]

\newtheorem{corollary}[theorem]{Corollary}

\newtheorem{lemma}[theorem]{Lemma}

\newcommand{\beq}{\begin{equation}}
\newcommand{\eeq}{\end{equation}}
\newcommand{\bcl}{\begin{center}}
\newcommand{\ecl}{\end{center}}
\newcommand{\re}{\mathbb{R}}

\newcommand{\iy}{{\infty}}

\newcommand{\ve}{\varepsilon}

\newcommand{\ms}{\medskip}

\newcommand{\pa}{\partial}
\newcommand{\Fn} {\rule{2mm}{2mm}}

\begin{document}

\null\thispagestyle{empty} \vspace{0.1in}
\begin{center}
{\Huge \textbf{Nonlinear Diffusion and}}\\ {\ }
\\{\Huge \textbf{ Image Contour Enhancement}} 
\end{center}
\begin{center}
\vspace*{0.5in} {\Large {\sc G.I. Barenblatt}\\
University of California at Berkeley, USA}
\end{center}
\begin{center}
{\Large \sc and J.L. V\'azquez}
\\{\Large Universidad Aut\'onoma de Madrid, Spain}\\
\end{center}

\begin{center}\vspace*{.2in}
\textsc{\large\today}
\end{center}

\begin{center}\vspace*{.5in}
\textsc{Abstract}
\end{center}

\noindent  The theory of degenerate parabolic equations of the
forms
\[
u_t=(\Phi(u_x))_{x} \quad \mbox{\rm and} \quad
v_{t}=(\Phi(v))_{xx}
\]
is  used to analyze the process of contour enhancement in image
processing, based on the evolution model of Sethian and Malladi.
The problem is studied in the framework of nonlinear diffusion
equations. It turns out that the standard initial-value problem
solved in this theory does not fit the present application since
it it does not produce image concentration. Due to the degenerate
character of the diffusivity at high gradient values, a new free
boundary problem with singular boundary data can be introduced,
and it can be solved by means of a non-trivial problem
transformation. The asymptotic convergence to a sharp contour is
established and rates calculated.

\vskip .3in



\vskip 1cm

\newpage


\section{\textsf{Introduction. A model for contour enhancement}}
\label{sect-intro}
\setcounter{equation}{0}

 This paper is devoted to study the behaviour of the solutions of
some degenerate parabolic equations with moving boundaries 
which appear in describing the technique of {\sl contour enhancement} in image processing.
Indeed, computer vision has become in recent decades a
mathematical discipline which relies on the differential-geometric
approach. More specifically, an appropriate technique of image
processing consists of formulating a partial differential equation
of evolution type for the {\sl image intensity}, $u(x,y)$. This
function, also called the {\sl grey level}, takes values in the
interval $0\le u\le 1$ and is defined on a two-dimensional  image
domain, $\Omega$. The usual evolution model leads to a nonlinear
equation of parabolic type, possibly degenerate or singular. The
nonlinearity is created by the law relating the image intensity
flux to the image intensity.

It has been observed by Perona and Malik \cite{PM}, 1990, that for a
suitable choice of this fundamental nonlinearity there appears an
effect of enhancement of image edges that has a strong interest in
the application to processing, denoising and recognition of images.
The Perona-Malik Anisotropic Diffusion model has had a deep influence
in the field, being at the source of many later developments. The model
proposed by Malladi and Sethian \cite{MS} leads (after
proper scaling) to the following equation for the image intensity:
\begin{equation}\label{eq1a}
u_t=(1+|Du|^2)^{1/2}\,\kappa
\end{equation}
where \ $D u$ \ denotes the spatial gradient of $u$, and $\kappa$
denotes the curvature of the surface $z=u(x,y)$. The equation
represents movement by curvature (curvature flow) and can be
written as
\begin{equation}\label{eq1b}
u_t=\frac{(1+u_y^2)u_{xx}- 2u_xu_yu_{xy}+(1+u_x^2)u_{yy}}{1+u_x^2
+ u_y^2}.
\end{equation}
We consider here the more general flow given by equation
$u_t=(1+|Du|^2)^{(1+p)/2}\kappa, $ where $p$ is a constant
parameter. In other words, we study the equations
\begin{equation}\label{eq2}
u_t=\frac{(1+u_y^2)u_{xx}- 2u_xu_yu_{xy}+(1+u_x^2)u_{yy}}{(1+u_x^2
+ u_y^2)^{1 +\alpha}}
\end{equation}
with parameter \ $\alpha=-p/2$. Along with the former case
$p=\alpha=0$, the case $p=-2$, $\alpha=1$ has also attracted the
attention of researchers (Beltrami flow, cf. Sochen et al.
\cite{SKM}).

The asymptotic treatment of these models done in the papers
\cite{MS, SKM} shows the enhancement of the intensity contrasts by
formation of regions of large intensity gradients, i.e., the
normal component of the image intensity gradient becomes quite
large. This phenomenon was analyzed in \cite{B2}, where  a further
simplification of the model was proposed in order to focus on the
{\sl boundary layer} where large gradients concentrate. Arguing
locally around a sharp gradient point and choosing the $x$-axis as
the direction normal to the boundary layer or front, we may
disregard the effect of $y$ derivatives with respect to the $x$
derivatives in (\ref{eq2}). In this way we get the {\sl reduced
equation}, which is just the one-dimensional version of
(\ref{eq2})
\begin{equation}\label{eq3}
u_t=\frac{u_{xx}}{(1+u_x^2)^{1 +\alpha}},
\end{equation}
where we have neglected $u_y, u_{yy}$. Different dimensional
constants appear in the model,  but they have been
scaled to unity here without loss of generality. The mathematical
problem consists in solving this equation with suitable boundary
data, namely, $u=0$ on the left-hand side of the contour and $u=1$
on the right-hand side (be that a finite or an infinite distance),
and initial conditions
$$ u(x,0)= u_0(x),
$$
 satisfying $0< u_0
<1$ and $u_0'>0$ in an interval $I=(a,b)$ and constant values
otherwise, zero to the left, 1 to the right. As was pointed out in
\cite{B2}, the phenomenon of {\sl gradient enhancement} takes place
in this model (in a proper setting) for all $\alpha \ge 0$: the spatial gradient of
the solutions, $u_x$, increases with time, and its support shrinks. 
Indeed, we  show below that such a  behaviour can be observed
in the larger exponent range  $\alpha>-1/2$. 
The conditions on $u_0$ can be relaxed, but then problems arise. Thus, the
less stringent size restriction $0\le u_0\le 1$ or the lack of
monotonicity create interesting alternatives, that will also be
briefly discussed.

Let us remark that this is not the only model which uses nonlinear
diffusion equations for image processing and enhancement. We refer
for alternatives to the works of Alvarez, Lions and Morel \cite{ALM}
and Caselles et al. \cite{Cas}, where further references are found. Let us
point out that in \cite{ALM} a fundamentally different assumption
is made, namely that the flux is perpendicular to the gradient,
and the basic equation is substantially different from the class
(\ref{eq2}) treated here. On the other hand, we mention that the
mathematical difficulties of the original model by Perona and
Malik have been further investigated by several authors, like  \cite{Kc} and
\cite{WB}. The application of Nonlinear Diffusion  to Image Processing is a
every active concern with many issues being discussed in the literature,
cf. e.g. \cite{KB, Weick}.


\section{\textsf{Asymptotic self-similarity and enhancement}}
\label{sect.as}
\setcounter{equation}{0}

Evolution equations like (\ref{eq3}) and many other variants have
been studied and are known in the literature under the general
name of nonlinear parabolic equations of diffusion type, or {\sl
nonlinear diffusion equations} for short. They are typically used
in describing processes of mass diffusion or thermal propagation.
A quite general one-dimensional form popular among PDE experts is
$$ u_t=(a(x,t,u,u_x))_x $$ with suitable conditions to make it
parabolic, at least in a formal sense, like $\partial a/\partial
u_x>0$. However, many practical applications (as in the present
case) involve functions $a$ for which $\partial a/\partial u_x\ge
0$, but the values zero or infinity can also be taken, and then
the equations are known as degenerate parabolic or singular
parabolic resp., cf. \cite{DiB, Kal}. In any case, a general
feature of this wide class of equations is the diffusive
character, which roughly means the spreading of the level sets of
the solutions as time advances. This property goes squarely
against the desired enhancement, therefore an extra mechanism must
be present if enhancement is to occur. We recall that in the Perona-Malik model
this mechanism was negative diffusion.

In the range of exponents $\alpha\ge 0$, equation (\ref{eq3}) falls
into the class of degenerate parabolic equations with degeneracy at $u_x=\infty$,
precisely the limit value which is of concern in contour enhancement.
An investigation of the phenomenon of gradient enhancement in
equation (\ref{eq3}) is performed in the paper \cite{B2} by
relating it to the convergence toward self-similar asymptotics of
an approximate model. It goes as follows: the observed evolution
of the solutions towards a configuration with large gradients
makes it plausible to further simplify the expression $1+u_x^2$ in
(\ref{eq3}) into $u_x^2$, so that the relevant reduced equation
becomes
\begin{equation}\label{eq4}
u_t=u_x^{-2(1+\alpha)}\,{u_{xx}}\,.
\end{equation}
It is further observed that  self-similar solutions for this
equation with end-levels $u=0$ and $u=1$ exist for all $\alpha>0$
and exhibit the similarity form
\begin{equation}
 u(x,t)=
F(\xi), \qquad \xi=(x-x_0)\,(t+t_0)^{1/2\alpha}\,,
\end{equation}
\noindent where $x_0$ and $t_0$ are parameters to be fixed, and
the profile $F(\xi)$ is an increasing function joining the levels
$F=0$ at a finite distance $\xi=-c<0$ to the level $F=1$ at
$\xi=c$. At these levels, which are taken at a finite distance
$\xi=\pm c$, the gradients are infinite (actually, the analysis in
\cite{B2} deals with {\sl decreasing} profiles joining $u=1$ to
$u=0$ but these two problems are obviously equivalent after a
mirror symmetry. We have chosen the increasing option to avoid
chasing around many minus signs).
 The scaling implies that
$$
u_x(x,t)= (t+t_0)^{1/2\alpha} F'(\xi), \qquad
\xi=(x-x_0)\,(t+t_0)^{1/2\alpha}\,,
$$
which shows that the solution is concentrated in an increasingly
narrower strip
\begin{equation}
S= \{(x,t): |x-x_0|\le c\,(t+t_0)^{-1/2\alpha}\}
\end{equation}
with gradients that diverge like $t^{1/2\alpha}$ as $t\to \infty$.
Let us remark that the asymptotic divergence of the gradients is
known to be exponential for $\alpha=0$, but the form is not
self-similar of the same type, as we will see below.
An important feature to be remembered of these solutions, along with the
 divergence as $t\to\infty$, is the {\sl infinite gradient condition} at the
 end-points of the domain of definition.
On the other hand, these self-similar solutions represent the
intermediate asymptotics of the problem with more general data.
This is demonstrated numerically in \cite{B2} when  infinite flux
data are imposed at moving end-points located at finite distance.

With these preliminaries we are ready to attack the construction
and analysis of solutions with steep fronts for a general class of
equations which includes  (\ref{eq4}) and  (\ref{eq3}) and is
natural for our application, and with the general data mentioned in
the Introduction.  We use the theory  of nonlinear diffusion equations
which has been strongly developed in the last decades.
It turns out that the standard initial value problem solved in
this theory does not suit our model, since it leads to dispersion
instead of concentration of the image.  A new {\sl free boundary
problem with singular boundary data} is then introduced and solved
by means of a number of non-trivial problem transformations
available for nonlinear diffusion problems, a subject that is developed
in detail in \cite{Vdarcy}. The concentration property of these solutions
is possible thanks  to the degenerate character of the diffusivity for
high gradient values. 

Our theory below covers existence and uniqueness of solutions of
suitable problems, existence and behaviour of the bounding
interfaces, and large-time behaviour. In particular, we obtain the
rates of convergence of the interfaces, an important question for
image enhancement.


\section{\textsf{Nonlinear diffusion equations}} \label{sectnde}
\setcounter{equation}{0}

As we have just said, we are going to construct a mathematical
theory for the one-dimensional evolution problem with initial and
boundary conditions of the type mentioned above, but replacing the
special class of equations (\ref{eq3}) by the  wider class of {\sl
nonlinear diffusion equations} of the form
\begin{equation}\label{eq5}
u_t=(\Phi(u_x))_{x}\,.
\end{equation}
The nonlinear function $\Phi$ expresses the dependence of  the
image flux on the gradient of the image intensity. It can be
called the {\sl flux law} of our process, or also the constitutive
function. In this general setting, $\Phi$ can be any increasing
real function defined in a suitable interval of gradients. We will
assume  for convenience that

\noindent {\bf (H${}_1$)} \quad  {\sl $\Phi(s)$ is $C^\infty$
smooth and strictly increasing in the intervals $0<s<\infty$ and
$-\infty<s<0$, allowing for a quite arbitrary behaviour as $s\to
0$ or $|s|\to \infty$.}

In  our opinion, this level of generality instead of obscuring the
problem makes more apparent the relation between nonlinear
diffusion and contour enhancement. Slightly more general
assumptions can be made on $\Phi$, but they are unessential for
our present purposes, cf. Section \ref{sect-dc}.

We want to  characterize the class of those $\Phi$ for which we
can construct  solutions with gradient-enhancement. We also want
to describe the rates at which the formation of steep profiles
takes place, thus justifying and extending the results of
\cite{B2} to the actual equations like (\ref{eq3}),  and to general
data. The mathematical problem we pose consists in solving the
nonlinear parabolic equation (\ref{eq5}) with initial conditions
\begin{equation}\label{eq.bc2}
u(x,0)= u_0(x)
\end{equation}
satisfying $0\le u_0\le 1$ and other suitable conditions (see
below), and boundary data $u=0$ on the left-hand  side of the
contour and $u=1$ on the right-hand side. Regarding the initial
data, we remark that we are mainly interested in monotone
solutions, i.e., $u_x\ge 0$, for equation (\ref{eq3}), and this
condition will naturally follow from a similar monotonicity
condition on the initial data. As for the boundary conditions, it
turns out that, depending on the form of  $\Phi$, the boundary can
be chosen to be located either at infinity or  at a finite
distance. This latter case will be the one of interest for us, and
then the problem must be properly posed as a {\sl free-boundary
problem}.

Let us recall the simplest examples of $\Phi$ and some of the
difficulties we will encounter. Indeed, a quite important and
simple example is the power function, that we
write as $\Phi(s)=(1/m) s^m$, defined for $s\ge 0$, so that it
agrees with (\ref{eq4}) with $m=-1-2\alpha$ since
$\Phi'(s)=s^{m-1}$. There is in principle no reason to restrict
the generality of the exponent $m$ in the mathematical treatment
to follow, and this will lead to quite different behaviour types
inside this family. The case $m=0$ is included in the form
$\Phi(s)=\log(s)$, i.e., $\Phi'(s)=1/s$. Note that in the  cases
with $m<0$ the function $\Phi$ is negative, but the important
quantity for the parabolic character of the equation, i.e.,
$\Phi'(u_x)$, is always positive.
Finally, we note that equation (\ref{eq3}) corresponds to
$\Phi'(s)=  (1+s^2)^{-(1+\alpha)}$; it degenerates as $s\to\infty$
in the prescribed range $\alpha\ge 0$, even if $\alpha>-1$, but the equation
is is perfectly parabolic in regions of bounded  $u_x$. 

\subsection{Non-monotone solutions}

\noindent The above statements are made on the assumption that the solutions are
monotone, $u_x\ge 0$, which is not unjustified in our problem
setting but restricts the mathematical generality. When
considering non-monotone solutions it is customary in the
nonlinear diffusion literature to extend the power nonlinearity to
arguments $s=u_x<0$ in the simplest symmetric way:
\begin{equation}\label{power}
\Phi(s)=\frac 1{m}\,|s|^{m-1}s\,,
\end{equation}
which gives   $\Phi'(s)=|s|^{m-1}$, always nonnegative. This definition
 poses no
problem when $m>0$ and solutions with changing sign exist
corresponding to initial data with the same property. However, we
are interested in exponents $m<0$ (so-called {\sl very fast
diffusion} in the literature) where  $\Phi'$ is singular at $s=0$,
and the whole function $\Phi$ is no more monotone. The difficulty
has been studied in \cite{RV2} leading to the consequence that
solutions with changing sign do not exist, even in the weak sense
when $m\le 0$. This is a consequence of the singularity of the equation
at the level $u_x=0$ and does not affect  equation (\ref{eq3}).
Due to this obstruction, we will concentrate here on problems
with monotone solutions. 


\section{\textsf{Second formulation as a nonlinear diffusion equation }}
\label{sectnde2} \setcounter{equation}{0}

\noindent If we formally differentiate equation (\ref{eq5}) with
respect to $x$ and put $v=u_x$, we obtain  the equation satisfied
by the image intensity gradient:
\begin{equation}
\label{eq6} v_{t}=\Phi(v)_{xx}\,,
\end{equation}
which is usually called the {\sl nonlinear filtration equation},
NLFE, and is the  most  standard  class of nonlinear diffusion
equations studied in the literature.  We can call it in this
context the ``differentiated equation". Inversely, we can recover the
$u$ formulation from a solution $v(x,t)$  for (\ref{eq6})
 by means of the rule
\begin{equation}
\label{inv}
u(x,t)= c + \int_{\Gamma} ( v\,dx + \Phi(v)_{x}\,dt)\,,
\end{equation}
integrated along any curve $\Gamma$ in the domain of definition of
$v$ which joins a fixed point $(x_0,t_0)$ to the generic point
$(x,y)$, cf. \cite{RV1, Vdarcy}. The constant $c$ is the value of $u$
at $(x_0,t_0)$ to be chosen at will in principle. The calculation
is justified for smooth solutions $v$ and smooth $\Phi$, but holds
in a much wider context.

Monotone solutions $u$ for equation (\ref{eq3}) translate into
nonnegative solutions $v$ for equation (\ref{eq6}), and conversely.
In this context, $u$ is usually viewed as the {\sl mass function}
for $v$, since when we take   $\Gamma$ to be a segment of the line $t=t_0$,
formula (\ref{inv}) becomes
$$
u(x,t_0)-u(x_0,t_0)=\int_{x_0}^x v(y,t_0)\,dy.
$$
The phenomenon of gradient enhancement can be then translated into
usual diffusion parlance as {\sl mass concentration}.
We will keep in the sequel the denomination {\sl intensity} or
{\sl image intensity} for the solution $u$ of equation
(\ref{eq5}), and we will view the solution $v=u_x$ of (\ref{eq6})
as the {\sl intensity gradient}. In that context
$\Phi(u_x)=\Phi(v)$ is the {\sl intensity flux.} Finally,
$\Phi'(u_x)=\Phi'(v)$ is the {\sl diffusivity.}


\section{\textsf{Solutions of Type I}} \label{sectt1}
\setcounter{equation}{0}

We now try to apply the more standard theory of nonlinear
diffusion equation to solve the evolution problem motivated by the
application in the Introduction. We show next that this
implementation can be performed in a rather standard way for a wide
class of $\Phi$'s, but such a process  gives an evolution with no concentration
effect, hence useless for our purposes in image processing.
The reader already familiar with diffusion theory may check the
contents of Theorem \ref{thm.1} and skip the rest of this section.

\subsection{Simplest setting}

\noindent Generalizing the well-known properties of the heat
equation (the choice $\Phi(s)=s$), we consider first the case
where  $\Phi$ is a $C^1$ function and $\Phi'$ does not vanish. We
pose the Cauchy problem for equation (\ref{eq5}) in the whole line
$x\in \re$ for $t>0$ with  bounded initial data $0\le u_0\le 1$,
and get a unique smooth solution $u(x,t)$ defined in
$Q=\re\times (0,\infty)$ and such that $0< u < 1$, $u$ is smooth
for all $x\in\re$ and $t\ge \tau>0$. We can also work with
equation (\ref{eq6}): then, if $v_0=u_{0,x}$ is locally integrable
(a Radon measure will also do) and such that
$$ \int_{\re}
v_0(x)\,dx=1,
$$
the solution $v(x,t)$ is smooth in $Q$ and
satisfies the same type of integrability condition, $ \int_{\re}
v(x,t)\,dx=1\,. $ For solutions $v\ge 0$ this is called
conservation of mass.  Using formula (\ref{inv}) we obtain a
smooth solution $u(x,t)$ in $Q$ such that the initial condition
$$
u(x,0)=u_0(x)=\int_{-\infty}^x v_0(y)\,dy
$$
holds, as well as the end conditions
$$
\lim_{x\to -\infty} u(x,t)=0, \quad \lim_{x\to
\infty} u(x,t)=1,
$$
hold locally uniformly in time. On the other
hand, $v_0\ge 0$ implies that $u_x=v\ge 0$,  which by the Strong
Maximum Principle implies $u_x>0$. It follows that $0<u(x,t)<1$ in
$Q$. Note that, due to the smoothness condition on $\Phi$, these
solutions are smooth in $(x,t)$ for all $x\in \re$ and $t>0$, even
if the initial data are not.

\subsection{Dispersion of solutions with time}

\noindent At first glance, the type of solution we have
constructed seems to solve our image problem. However, it lacks a
basic ingredient,  i.e., the eventual concentration of intensity
gradients. On the contrary, if we consider the heat equation
$u_t=u_{xx}$, the solutions, which  can be expressed in terms of
error functions, spread in time and its gradients go to zero.
 for instance, if $v_0=u_{0,x}$ is integrable, then $v(x,t)$ goes
to zero as $t\to\infty$ at a rate of the order of $t^{-1/2}$, and
takes on a Gaussian space shape,
$$
v(x,t)t^{1/2}\sim \mbox{\rm exp} (-x^2/4t).
$$
And a similar result (with possibly different rates)
applies to more general data $u_0$ having definite limits at $\pm \infty$,
and to all the functions $\Phi$ of the above class.

\subsection{Extension of Type I solutions to other $\Phi$}

\noindent Even if this is not the class of solutions we are looking for, we
pursue a bit further the analysis, since it is the standard class
found in studies of nonlinear diffusion or thermal propagation.
The study will also serve for comparison with the ``correct"
solutions of Type II. However, the reader may choose to skip this
the rest of the section and proceed with Section \ref{sect-ct}.

The class of constitutive functions $\Phi$ for which
there exists a class solutions $0\le u\le 1$ with a dispersive
character can be extended to include a quite general choice of
$\Phi$ on the condition of allowing for a suitably generalized
concept of solution. Thus, it is well-known following B\'enilan,
Crandall and others \cite{Be, BeC}  that we may take as $\Phi$ any
continuous nondecreasing function and the Cauchy problem is then
well-posed in the class of so-called {\sl mild} solutions with
$L^1(\re)$ data. Actually, $\Phi$ can be allowed to be
discontinuous but this will be of no particular interest here; on
the contrary, we will keep the assumption of smoothness
and strict monotonicity for $0<s<\infty$ for simplicity of
presentation (and since it is satisfied in the application we
are dealing with).


Let us comment on  the main properties of the solutions for this
class of equations. For the case of power nonlinearities $\Phi(s)=s^m/m,$ $s\ge 0$,
mentioned above, rather complete details are known,  \cite{Ar, Kal, Vb}. We can
consider that $m$  is a real parameter, positive in principle but
not necessarily as we will see. The definition is  extended to
solutions with negative values as in formula (\ref{power}). Thus,
for $m>0$ equation $v_t=(v^m)_{xx}$ (or better,
$v_t=(|v|^{m-1}v)_{xx}$) generates a positive semigroup of
contractions in the space $L^1(\re)$. In other words, for every
$v_0\in L^1(\re)$ there exists a unique function $v\in
C([0,\infty:L^1(\re))$, \ $v\ge 0$ with $v^m\in L_{loc}^1(Q)$ such
that the equation is satisfied in the sense of distributions in
$Q=\re\times (0,\infty)$, the initial data are taken in the
$L^1(\re)$-sense, and  the map $v_0\mapsto v(\cdot,t)$ is an
$L^1$-contraction. Moreover, the total mass is conserved $$
\int_0^{\infty} v(x,t)\,dx= \int_0^{\infty} v_0(x)\,dx. $$ In our
application we still have to impose the extra condition of total
mass 1, and then we recover the image intensity by means of the
formula $$ u(x,t)=\int_{-\infty}^{x} v(x',t)\, dx', $$ and the
intensity level goes from $u=0$ at minus infinity to $u=1$ at
infinity. The asymptotic behaviour of these solutions has been
carefully calculated in the literature, cf. e.g. \cite{Ang2, FrKa80,
Vncp}: solutions $v(\cdot,t)$ with finite mass go zero uniformly as
\begin{equation}\label{as.1}
v(x,t)\sim t^{-\gamma}, \quad \gamma=\frac1{m+1},
\end{equation}
which in the notation of the Introduction
means $\gamma=-1/2\alpha$. Actually,  the asymptotic rate comes
from comparison with the source-type self-similar solutions
\cite{Barbook} which take the form
$$
v(x,t)=t^{-\gamma}G(x\,t^{-\gamma})
$$
for a certain symmetric profile function $G\in C_b(\re)$ such that
$G(\xi) \to 0$ as $\xi\to \infty$. This formula suggests that the
behaviour will be maintained as long as $\gamma>0$, hence, as long
as $m>-1$, and this turns out to be true (for the study when
$-1<m\le 0$ \ cf. \cite{ERV}). Indeed, this exponent range is
optimal since there are {\sl no solutions} with finite integral
for $m\le -1$.

Let us look a bit closer at the kind of initial data that we want
to consider in the outmost generality. We recall that we want
$u_0$ to be increasing and bounded between 0 and 1. Hence,
$v_0=u_0'$ has to be defined and locally integrable, or at most be
a measure, in an interval $(a,b)$ with possible divergence at the
end-points, but anyway with finite integral. Existence of
solutions (with $v_0$ extended all of $\re$ with value 0 otherwise) offers no problem in
all the range  $m>-1$ and the solutions are bounded for all $t>0$.
But large values will be the origin of the new class of solutions
to discuss in the following sections.

By integration we obtain a solution of the problem
$$
 ({\rm P}_{I}) \qquad
 \left\{\begin{array}{ll}
 u_t=\Phi(u_x)_{x} \qquad & \mbox{ in } \ Q=\re\times (0,\infty)\\
 u(x,0)=u_0(x)  &\mbox{ for } \  x\in \re\\
 u(x,t)\ge 0,
 \end{array}\right.
$$
where $u_0$ is any nondecreasing continuous real function joining
the levels $u_0=0$ at $x=-\infty$ to $u_0=1$ at $x=\infty$. This
is what we term as Solutions of Type I. Note that integration of
the self-similar profile gives for $u$ the form
$$
u(x,t)= F(x\,t^{-\gamma})
$$
where $F$ is a primitive of $G$ with $F(-\infty)=0$,
$F(\infty)=M$ that we want to normalize to $M=1$.

\begin{theorem}\label{thm.1} Let us consider the initial-value
for equation {\rm (\ref{eq5})} posed in $Q=\re\times (0,\infty)$
 with power function $\Phi(s)=s^{m}/m$. If $m>-1$ then for every
 nondecreasing $u_0$ with $u_0(-\infty)=0$,  $u_0(\infty)=1$,
there exists a unique continuous weak solution in the sense
$u(x,t)\ge 0$ such that $u(\cdot,t)$ jumps from $0$ to $1$ as $x$
ranges over the line $x\in\re$. If $m\le 0$ the last
condition is essential to ensure uniqueness. This class  of solutions
has bounded gradients for strictly positive times ($t\ge \tau>0$),
 and spreads in space as time
advances and the asymptotic decay formula {\rm (\ref{as.1})} holds for $v=u_x$.

On the contrary, if $m\le -1$ solutions for this initial-value
problem with bounded data do not exist.
\end{theorem}

The results of the power case reflected in this theorem  can be
generalized to equations with  more general $\Phi$. Only the
behaviour of $\Phi(s)$ at $s=0$ and $s=\infty$ will determine the
different behaviour types.  In order to tackle the first, we assume
that the  initial data (and  hence the solutions) $v$ are bounded.
Then it is known that the  condition for existence with finite
mass is
\begin{equation}\label{}
 \int_0^s \Phi'(s)s\,ds\le \infty\,.
\end{equation}
The question of large arguments is similar to the power case. These issues will be discussed at length
in \cite{Vbk}. For a study of these questions in the class of self-similar solutions
cf. \cite{FV}.

\subsection{Existence of sharp interfaces}

\noindent There is an interesting subclass of equations (\ref{eq5}), or
(\ref{eq6}),  where sharp interfaces appear. Let us look first at power
 nonlinearities. Indeed, for exponent $m>1$ solutions of the Cauchy problem
for equation (\ref{eq5}) with initial data having compact support
will keep this property for all times, while infinite propagation
occurs whenever $m\le 1$ cf. \cite{BV, Kal}. In the first case,
given an integrable function $u_0\ge 0$ with integral 1 the
solution of the problem can be seen as the a classical solution
equation (\ref{eq5}) in a domain $$ \Omega= \{(x,t):
-l(t)<x<r(t)\}\,, $$ with initial conditions $u(x,0)=u_0(x)$ and
boundary conditions $$ \left\{\begin{array}{l} u(x,t)=0, \quad
\Phi(u_x)=0 \qquad \mbox{ for } \ x=l(t),\\ u(x,t)=1, \quad
\Phi(u_x)=0 \qquad \mbox{ for } \ x=r(t).
\end{array}\right.
$$
The lines $x=l(t)$ and $x=r(t)$ are called interfaces or moving
boundaries and are completely determined by the over-specified
conditions. They are known to be smooth (analytic) functions of
$t$ \cite{Ang, AV}, and diverge as $t\to \infty$ like $O(t^{\gamma})$,
thus giving a quantitative estimate of the dispersion effect. We remark in
passing that the presence of interfaces means also that the
equation is not uniformly parabolic at those points, and consequently
the solutions have limited regularity.

On the other hand, when $m\le 1$ the same Cauchy problem leads to
positive solutions with  $l(t)=-\infty$ and $r(t)=\infty$. The property of
null flux is equivalent to imposing $u_x=0$ at $\pm \infty$, a quite
natural condition in view of the values $u=0,1$ at $\pm \infty$.
This condition is automatic for $m>0$. However, for $m\le 0$ we can
have {\sl new solutions with decreasing total mass}, i.e.,
such that
$$
\frac{d}{dt}\int v(x,t)\,dx <0,
$$
and they can be even forced to vanish identically in finite time by
controlling the outgoing flux at $x=\pm\infty$. We refer to \cite{RV93,RV1}
for a detailed analysis. In the integrated version
they would lead to solutions with a restricted grey range.
For the more general class of $\Phi$ mentioned above the condition
to have finite interfaces is
\begin{equation}\label{}
 \int_0^s \frac{\Phi'(s)}{s}\,ds\le \infty.
\end{equation}
The fact that this condition is necessary and
sufficient can easily seen on the family of traveling waves. The
property was first pointed out in \cite{OKC}.


\section{\textsf{ Solutions of Type II. Conjugate formulations}}
\label{sect-ct}
\setcounter{equation}{0}

We now address the main question of this paper, the construction of solutions with large
gradients, appropriate for the contour  enhancement problem. 

\subsection{ Basic Free Boundary Problem. \rm\normalsize It is formulated as
follows:}
\vskip -.3cm

{\sl Given an increasing function
$u_0(x)$ defined in an interval $(a,b)$ with end values $u(a+)=0$,
$u(b-)=1$, to find a continuous function $u(x,t)$ and continuous
curves $x=l(t)$ and $x=r(t)$ such that

 {\rm (i)} $l(0)=a$, $r(0)=b$, and $l(t)<r(t)$ for some
time interval $t\in (0,T)$.

{\rm (ii)} $u$ solves the following problem in $\Omega=\{(x,t):
0<t<T, \ l(t)<x<r(t)\}$:
$$
 ({\rm P}_{\rm II}) \qquad \qquad
 \left\{\begin{array}{ll}
 u_t=\Phi(u_x)_{x} \qquad & \mbox{ in } \ \Omega\\
 u(x,0)=u_0(x)  &\mbox{ for } \ a\le x\le b\\
 u(l(t),t)=0, u_x(l(t),t)= +\infty &  \mbox{ for } \ 0<t<T,\\
 u(r(t),t)=1, u_x(r(t),t)= +\infty &  \mbox{ for } \ 0<t<T.
 \end{array}\right.
$$}
Such a triple $(u,l,r)$ is called a solution of Type II. The
regularity required from $u$ as a solution of the problem will
depend on the generality of the data. At least $u$ will be
continuous in the closure of $\Omega$. Furthermore, to avoid
unnecessary generality we will ask $u$ to be smooth
in the interior of $\Omega$. Finally, the requirement of
monotonicity is not intrinsic from the mathematical point of view,
but it suits the application and allows for the use of the
powerful conjugate formulations. Here is the existence result that
we are going to prove:

\begin{theorem}\label{thm.tII} Let  $\Phi$ be a flux function that is
defined, smooth and with $\Phi'(s)>0$ for all $s>0$. Assume moreover that
 $\Phi(\infty)$ is finite. Then for every increasing function $u_0(x)$ defined in an
interval $[a,b]$ with $u(a)=0$, $u(b)=1$ and $u_{0}'\ge c>0$ there
exists a unique continuous function  $u(x,t)$  which is defined in
a set $\Omega$ as above, is smooth  and strictly monotone in $x$
for $0<u<1$, and there exist continuous curves $l(t)$ and $r(t)$,
such that the triple {\rm ($u,l,r$)} solves  problem ${\rm
(P}_{\rm II}{\rm )}$ in $\Omega_T$. Besides,  $u_x\ge c >0$
whenever $0<u<1$. \end{theorem}

We will also show in our construction  that the boundary curves
are monotone: $l(t)$  is nondecreasing,  $r(t)$ nonincreasing. On the other hand, there
is the problem of determining whether $T$ if finite or infinite. This depends on
$\Phi$ as we will see below. 

\subsection{Conjugate formulation}

\noindent When dealing with smooth monotone solutions $u_x>0$ we can invert
the variables $x$ and $u$ and
write $x=X(u,t)$. Then $u_x\cdot x_u=1$, and after some
computations we get the partial differential equation satisfied
by $x$ a a function of $u$ and $t$:
\begin{equation}
\label{eq7}
x_{t}=(\Psi (x_u))_{u}\,,
\end{equation}
where $\Psi$ is the {\sl conjugate flux function} (conjugate to
$\Phi$), defined for $s>0$ as
\begin{equation}
\label{eq.conj}
\Psi(s)=-\Phi(1/s)\,.
\end{equation}
Differentiation of equation (\ref{eq7}) with respect to  $u$ gives rise to  the {\sl differentiated conjugate equation} for
$w={\partial X}/{\partial u}$  as a function of $u$ and $t$: 
\begin{equation}
\label{eq8}
w_t =(\Psi (w))_{uu}\,.
\end{equation}
We complete the list of related equations with the direct
differentiated equation for $v=\partial u/\partial x$, already
seen:
\begin{equation}
\label{eq8b}
v_t =(\Phi (v))_{xx}\,,
\end{equation}
and then $ v= 1/ w$. It is important to point out that these relations are equivalent
to the well-known  B\"acklund transform \cite{BK} between the main variables
$v$ and $w$ of the second and fourth formulations. Indeed, we have
$$
v(x,t)=\frac1{w(u,t)}, \quad u(x,t)=c + \int_{\Gamma} ( v\,dx +
\Phi(v)_{x}\,dt)\,.
$$
The reader is referred to \cite{RV1, Vdarcy, Vposit} for other applications of this
technique.

\subsection{Posing and solving the conjugate problems}

\noindent We now show how to use the conjugate formulations to
solve the original problem. We assume that $\Phi$ a flux function defined for all
$s>0$ and such that $\Phi(\infty)$ is finite, say,
$\Phi(\infty)=0$. Then $\Psi(s)$ is continuous at $s=0$ and
$\Psi(0)=0$. This is the class of flux functions for which the
conjugate problem looks simpler. Since we have assumed that $\Phi$ is smooth,
so is $\Psi$ in its domain of definition.

\noindent (i)  Since we are interested in solving the conjugate
problem as an auxiliary step for Problem ${\rm (P}_{\rm II}{\rm )}$ we
will relate the initial values for the function $w(u,t)$ to $u_0$
as follows. Assuming that $u_0$ is continuous and strictly
monotone in the interval $I=\{a<x<b\}$, with $u_0(a)=0$,
$u_0(b)=1$, and $C^1$ smooth inside $I$ with $du_0/dx$ bounded
below away from zero, we define the inverse function
$x=h(u)=u_0^{-1}(u):[0,1]\to [a,b]$.
\begin{equation}\label{conj.in.dat}
w_0(u)=\frac{1}{u_{0,x}(h(u))}\,,
\end{equation}
which is defined for  $0\le u\le 1$ and is positive, bounded and
smooth inside, i.e., for $u\in (0,1)$.

\noindent (ii) We then solve the conjugate problem
$$({\rm P}_{\rm c}) \qquad \qquad
 \left\{\begin{array}{ll}
 w_t=\Psi(w)_{uu} \qquad & \mbox{ for } \ 0<u<1, t>0\\
 w(u,0)=w_0(u)  &\mbox{ for } \ 0\le u\le 1\\
 w(u,t)=0 &  \mbox{ for } \ u=0,1.
 \end{array}\right.
$$
As initial data we choose a nonnegative, bounded function $w_0$.
Under these conditions Problem ${\rm (P}_c{\rm)}$ has a unique
solution by virtue of well-known  nonlinear parabolic theories described
in the previous section, cf. \cite{Be, BeC, BeG}; but note that now we are
dealing with the homogeneous Dirichlet problem. The solution can be obtained as limit of
the solutions $w_\ve(u,t)\ge \ve$ of the nondegenerate problems
with initial data $w_{0,\ve}(u)=w_0(u)+\ve$, $\ve>0$. In the
monotone limit we get
$$ \lim_{\ve \to 0} w_\ve(u,t)=w(x,t)
$$
which is nonnegative, continuous and bounded. Under the additional
assumption that $w_0$ is locally bounded away from zero, it is
easily proved that the solution $w(u,t)$ is positive, hence
classical, in a strip
$$
S_T=\{(u,t): \ 0<u<1, \ 0<t<T \}.
$$

\noindent (iii)  We also need to know something about the
large time behaviour of the solutions to this problem. In the full
generality ($\Psi$ continuous at $s=0$), it can be proved rather
easily that $w(x,t)$ goes to zero in uniform norm as $t\to
\infty$. However, the rate depends on $\Psi$, as we will see in
detail in the next sections. Indeed,  depending on $\Psi$ it may
happen that the solution vanishes identically
after a finite time $T>0$ (so-called {\sl finite-time
extinction}). For power nonlinearities $\Psi'(s)=s^{q-1}$ this
happens iff $ 0< q< 1$, \cite{BH}. On the other hand, it is well-known
that for $q\ge 1$ the decay rate is
$O(t^{-1/(q-1)})$, while for $q=1$ it is exponential.

\noindent (iv) Next, we pass to the integrated version using
the  formula
\begin{equation}\label{zboundary}
z(u,t)=\int_{\Gamma} w\,du + \Psi(w)_u\,dt,
\end{equation}
where $\Gamma$ is any piece-wise smooth curve in $(u,t)$ space
starting from a fixed point, say $u=1/2, t=0$ and arriving at a
generic point $(u,t)$. In this way we obtain a solution of the
integrated equation $z_t=(\Phi(z_u))_u $, much as we did in the
case of the original pair of formulations. Note that
\begin{equation}\label{zto0}
z(1,t)-z(0,t)=\int w\,du\to 0
\end{equation}
as $t\to T$.

\noindent (v) Let us examine the curves $z_{(u)}(t)=z(u,t)$ for
fixed $u\in [0,1]$. It
is clear from the smoothness of the solutions that these curves
are $C^\infty$ smooth for every $0<u<1$. We are interested in the
limit curves
\begin{equation}
l(t)=\lim_{u\to 0} z(u,t), \qquad
r(t)= \lim_{u\to 1} z(u,t).
\end{equation}
The limit is well-defined for every $0<t<T$ by monotonicity of $z$
as a function of $u$.  These curves will show up in the next subsection
as the interfaces of the original problem.

\medskip

\begin{lemma} \label{lemmaA}  {\sl The curves $z=l(t)$ and $z=r(t)$ are
continuous and monotone for all $0<t<T$ with $a\le l(t)<r(t)\le b
$ ($r(t)$ is nonincreasing and $l(t)$ nonincreasing). As $t\to T$
we have $r(t)-l(t)\to 0.$}
\end{lemma}

\noindent {\sl Proof.} For $0<s<t<T$ and $0<u_1<u_2<1$ we have
$$
z(u_1,t)-z(u_1,s)=\int_{u_1}^{u_2} w(u,t)\,du +
\int_{u_1}^{u_2} w(u,s)\,du + (z(u_2,t)-z(u_2,s)).
$$
Fix now $s>0$ and let $t$ be a bit larger than $s$. Since $w$ is
bounded uniformly for $t\ge s$ the two integrals are uniformly
small as long as $u_1-u_2$ is small. We will fix now $u_2\sim 0$
and let $u_1\in (0,u_2)$ go to 0. It is clear then that
$z(u_1,t)-z(u_1,s)$ is uniformly small and goes to 0 as $t\to s$.
Hence, $l(t)$ is continuous at $t=s$. The argument for $r(t)$ as
$u\to 1$ is the same.

For classical solutions the monotonicity of the limit curves is a
consequence of the estimate
$$
\frac{\partial  z}{\partial t}=\frac{\partial  }{\partial u} \Phi
(w)
$$
which must be positive at $u=0$ and negative at $u=1$ because
$w\ge 0$ and we have zero boundary conditions. For the general
case we use the dependence of the solutions on $\Phi$ as
demonstrated in \cite{BeC}. Hence, $l'(t)\ge 0$, $r'(t)\le 0$.
We will see another proof below.
Note that the rest of the curves $z_{(u)}(t)$ need not be monotone
unless $w_0$ is a rearranged symmetrical function.

As for the last statement, it follows immediately from \ref{zto0}.
It will mean for the original problem that the solution
concentrates into a vertical profile, thus showing the formation
of the desired front. \qquad \Fn

\subsection{Inversion. Solutions of Type II}

\noindent   Thanks to the fact that $\partial z/\partial u=w> 0$,
we can invert the dependence between $z$ and $u$ in the previous
construction to get a function $u=u(z,t)$ that is easily shown to
satisfy the equation
$$
u_t=(\Phi(u_z))_z \, .
$$
Besides, $u$ is a monotone function of $z$ and takes the values
$u=0$ and $u=1$ respectively at the left and right endpoints of
the domain of definition
$$
\Omega_z=\{ l(t)<z<r(t)\}, \quad l(t)= z(0,t), \quad r(t)=z(1,t).
$$
where $z(\cdot,t)$ is the function defined in (\ref{zboundary}).
Therefore, $u(z,t)$ is a candidate to solve our original problem if
we identify the independent variable $z$ with $x-c$, where $c$ is
uniquely determined by the relation \ $ u_0(c)=1/2.$

In order to check that we have solved the original problem we
still have to check some particulars.  It is clear that
$v=u_x$ is related to the original $w$ by the formula
$$
v(x,t)=\frac1 {w(u,t)}\,,
$$
which simply states the derivative rule for the inverse function,
and $u$ in the second member is given by $u(z,t)$, $z=x-c$, as
explained before. Here comes an important point: since $w$ takes
on zero boundary values, {\sl  $v(x,t)$ diverges at the endpoints
of its domain of definition,} $\Omega$. In other words, the
solutions of the original problem $u=u(x,t)$ enjoy the property of
infinite gradients at the endpoints of the strip where they are
defined. Since $\Phi(\infty)=0$ this also means zero flux at these
points, a reasonable requirement, which explains why this condition
has to be imposed on $\Phi$.

As for the initial data, we have the  mass formula
$$
x=\int_{1/2}^u w_0(u)\,du +  c, \qquad \mbox{with} \ u_0(c)=\frac12,
$$
so that $x$ ranges over an interval $[a,b]$ when $u$ goes from 0
to 1, i.e., $a=z(0,0)+c$, \ $b=z(1,0)+c$. This rule is accompanied
by the rule $ u=\int_{l(t)}^x v(x,t)\,dx.$

\subsection{Uniqueness} Uniqueness of our class of monotone solutions
works by translating any couple of
solutions of Problem ${\rm (P}_{\rm II}{\rm )}$ with the same
initial data into the conjugate formulation. They continue to have
the same initial data. Uniqueness of weak solutions is well-known
for that equation. The Theorem is proved. \Fn

\ms

\subsection{ Front formation} The asymptotic formation of a
vertical front is a simple consequence of the fact that there
exists the limit
$$
\lim_{t\to T} r(t)=\lim_{t\to T} l(t)= x_\infty\in (a,b).
$$
The existence of the common limit follows from Lemma \ref{lemmaA}.


\section{\textsf{Uniqueness and comparison for Type II solutions}}
\label{sect-uni}
\setcounter{equation}{0}

The previous construction provides for existence and uniqueness of monotone solutions.
We tackle next the property of comparison, which is stronger that
uniqueness.

\begin{theorem}\label{thm.tIIb}
Comparison applies to the solutions of Theorem {\rm \ref{thm.tII}} : let $u_1(x,t)$, $i=1,2$, be two solutions
defined in a strip $S_T=\re\times (0,T)$ having  initial data $u_{0i}$,
where the solutions and the data have been extended by 0 the the left of the
definition domain, by 1 to the right. If $u_{01}(x)\ge u_{02}(x)$, then
for all $x\in\re$ we have
\begin{equation}
u_1(x,t)\ge u_2(x,t) \quad \mbox{in } \ S_T\,.
\end{equation}
\end{theorem}

\noindent {\sl Proof.}   The main idea is to use a  contradiction
argument at  touching points of ordered solutions. We proceed in several stages
and  need a definition, taken from the theory
of viscosity solutions in \cite{CV}. We say that two solutions are {\sl strictly
separated} at a time $t$ if $u_1>u_2$ at all points intermediate to the
bounding interfaces, and these are also separated.

 (1) Given two solutions $u_1$ and $u_2$ with ordered data we must obtain approximations
to which the contradiction argument can be applied. First, we separate the initial data
of the solutions by displacing the first solution $u_1$ to the left (the problem is invariant under space
displacements). Next, assume that we have modified $u_1$ so that
it has a finite slope on the left-hand free boundary, and $u_2$ so that
it has a finite slope on the hand-side free boundary.  In this situation we argue on
the first point and time $(x_0,t_0)$
where the graphs of both solutions touch and discover that: (i) it cannot happen
with value $u\in (0,1)$ by
the Strong Maximum Principle; (ii) not on the left-hand side where $u_1=u_2=0$,
 because at such a point we have the right derivative: $u_{2,x}=+\infty$
while $u_{1,x}$ is finite, which contradicts the fact that we still must have
$u_1(\cdot,t_0)\ge u_2(\cdot,t_0)$; (iii) same argument applies to the right-hand
side where $u_1=u_2=1$. We conclude that such solutions $u_1$ and $u_2$ cannot
touch, hence strict separation is preserved in time, and
comparison $u_1\ge u_2$ holds in the strong sense, i.e., with strict separation
at all times.

 (2) We must now prove that the original solutions can be approximated by solutions
with finite derivatives on the lateral boundaries. This can be obtained by an easy
modification of the previous construction.  Thus, if $u$ is the
constructed solution, it can  be approximated at the level of the conjugate problem
 by putting the value $w=\ve $ on the left-hand side boundary $u=0$.
It is rather standard monotonicity argument that $w_\ve\to w$ in the limit $\ve\to 0.$
Undoing the transformation  this means that $u_\ve$ converges to the original
constructed solution $u$. A similar construction applies to $u_2$ by modifying
$w_2$ at the border $u=1$. After displacement of $u_\ve$ we get
$u_{1,\ve}(x+\delta,t)\ge u_{2,\ve}(x,t)$. Let now $\ve\to 0$. \quad \Fn

We derive next another proof of the monotonicity of the
interfaces.

\begin{corollary} The interfaces are monotone in all cases.
\end{corollary}

\noindent {\sl Proof.} Given a solution $u$ with initial interfaces
$a$ and $b$, we may place a very steep solution $\hat u$ to the
left, i.e., with $\hat u=0$ for $x\le a-2\ve$, $\hat u = 1$ for
$x\ge a-\ve$. If $\hat u'$ is symmetric then the interfaces are monotone,
hence the left interface of $\hat u$ lies to the right of $a-2\ve$.
Then $\hat u\ge u$ for all times, hence $l(t)\ge a-2\ve$, and in
the limit $l(t)\ge a$. The same argument applies taking origin of
times at any $y\in (0,T)$, hence $l(t)$ is monotonically
increasing in $(0,T)$. Likewise, $r(t)$ is proved to be monotone
decreasing. \quad \Fn


\section{\textsf{Self-similarity and asymptotics for power nonlinearities}}
\label{sect-ae}
\setcounter{equation}{0}

In the power case, where $\Phi(s)=(1/m) s^m$, the asymptotic
condition ``$\Phi(\infty)$ finite" means $m<0$. Then $\Psi(s)=(1/q)
s^q$ with $q=-m>0$. Reminder: we consider only monotone solutions of (\ref{eq5})
and nonnegative solutions of (\ref{eq6}), and $q=1 + 2\alpha$
in  the notation of the first section. Therefore, $\alpha>-1/2$
means $q>0$. 

There is a unique solution of problem ${\rm (P}_c{\rm )}$ with
initial data as in Section \ref{sect-ct}. The large time behaviour
depends on the exponent. For $q\ge 1$ solutions are defined and
nontrivial globally in time, while for $0<q<1$ they exist only in
a finite interval $0<t<T$ and $w(x,t)\to 0$ uniformly as $t\to T$.

As for the exact behaviour, if $\alpha\ne 0$ the self-similar behaviour
reduces to writing $U=\Phi(\xi)$ with $\xi=x\,t^{1/2\alpha}$ and integrating the ODE
 \begin{equation}
\Phi''= \frac1{2\alpha}\xi\,(\Phi')^{2+q}
 \end{equation}
 to obtain in the monotone region
  \begin{equation}
\Phi'(\xi)= A\,(K^2-\xi^2)^{-1/(2+2\alpha)}, \quad
A=\left(2\alpha /(1+\alpha)\right)^{1/(2+2\alpha)}\,.
 \end{equation}
 and a constant $K>0$ determined by the end values $\min\Phi(\xi)=0$,
 $\max\Phi(\xi)=1$.  For the large-time behaviour of general
solutions we have to make a study in different cases. We start
with $\alpha=0$ which was not included before.

\noindent {\bf Case $\alpha=0$, hence $q=1$.}  The conjugate equation is
the linear heat equation. It is well-known that the solutions  behave for
large times like
 \begin{equation}\label{w1a}
 w(u,t)\sim C_1\,W(u,t)\equiv C_1\,e^{-\lambda t}f_1(u),
\end{equation}
 where $\lambda=\pi^2$ is the first eigenvalue of the Laplacian in
 $[0,1]$, $f_1(x)=\sin(\pi u)$ is the first eigenfunction,
 and $C_1>0$ is a constant that depends on the
 initial data. Here and in the sequel we denote by capital letters the
quantities corresponding to self-similar solutions.
We  get the intermediate step
 \begin{equation}\label{w1b}
 x-c= z(u,t) \sim  Z(u,t)\equiv - C\,e^{-\lambda t}
 \,\cos(\pi u)\,,
\end{equation}
with $C=C_1/\pi$. This formula can be viewed as an implicit expression of \
$u=u(x,t)$ \ defined in the space between the interfaces. Note that $u=1/2$ for
$x=c$ and all $t>0$. Putting  $c=0$ we get for the self-similar
solution
 \begin{equation}\label{w1b2}
U(x,t) = F(x\,e^{\pi^2 t}/C)\,,
\end{equation}
where $F'$ behaves like $d^{-1/2}$ near the end-values or interfaces, $d$ being
the distance to
these points. Moreover, the interfaces
\begin{equation}\label{w1d}
R(t)=C\,R_1(t), \qquad L(t)=-C\,R_1(t), \qquad \mbox{\rm with} \
R_1(t)= e^{-\pi^2 t}\,.
\end{equation}
These estimates become the first order approximation when we consider
general solutions:
\begin{equation}
u(x,t)\sim F(x\,e^{-\pi^2 t}/C),\qquad
r(t), \ l(t)\sim c \pm C\,R_1(t).
\end{equation}
The gradients are given by
\begin{equation}\label{w1c}
 u_x(x,t)= \frac 1w \sim \frac1{C}\,e^{\lambda t}\frac1{\sin(\pi u)}\,.
\end{equation}
They  blow up exponentially as $t\to\infty$ and like $O(d^{-1/2})$
at the interfaces, where $d$ is distance to the interface.

\medskip

\noindent {\bf Case $\alpha>0$, hence $q>1$.} This is the case
treated in \cite{B2} by the direct self-similar method. We recover
the behaviour of general solutions from the conjugate problem. We
have a self-similar asymptotic expression in the form of
separation of variables:
 \begin{equation}\label{w2a}
 w(u,t)\sim W(u,t)\equiv t^{-1/(q-1)}f_q(u),
\end{equation}
 where $f_q(u)>0$ is the unique solution of the associated elliptic problem
$$
(f^q)'' + \mu f=0 \quad \mbox{in } \ [0,1], \quad \mu=
\frac{q}{q-1}\,,
$$
and $f=0$ at the endpoints. This problem was studied in \cite{AP}. In
this case there is no free constant, i.e., the behaviour is
universal. From it we get the intermediate asymptotic estimates
 \begin{equation}\label{w2b}
 x-c= z(u,t) \sim   t^{-1/(q-1)} g_q(u)\,.
\end{equation}
This formula defines implicitly $u=u(x,t)$ in the space  between
the interfaces as
\begin{equation}
u\sim F_q((x-c)\,t^{1/(q-1)}).
 \end{equation}
Since $f(u)$ behaves exactly like $O(u^{1/q})$ near the end point $u=0$, and
$O((1-u)^{1/q})$ near $u=1$, it follows after integration and inversion that
$u=F$ behaves like $O(d^{(q/(q+1)})$,
where $d$ is the distance to the interfaces. These are given by
\begin{equation}\label{w2d}
l(t)=c- C_q \,t^{-1/(q-1)} , \qquad r(t)=c+  C_q\, t^{-1/(q-1)}.
\end{equation}
The gradient is given by
\begin{equation}\label{w2c} u_x(x,t)\sim t^{1/(q-1)}\frac1{f_q(u)},
\end{equation}
and blows up like $O(t^{1/2\alpha})$ as $t\to\infty$
and like $O(d^{-\gamma})$ at the interfaces, where $\gamma=
1/(q+1)= 1/(2\alpha+2)$.  

\medskip

\noindent {\bf Case $-1/2<\alpha<0$, hence $0<q<1$.} This
is a new case, {\sl not included} in the modelization of \cite{MS, B2}.
We get solutions from the conjugate problem which exist for a
finite time $T$ and behave as $t\to T$ like the separation of variables
formula:
 \begin{equation}\label{w3a}
 w(u,t)\sim (T-t)^{1/(1-q)}f_q(u),
\end{equation}
 where $f_q(x)>0$ is the unique solution of the associated elliptic problem
$$
 (f^q)'' + \mu f=0 \quad \mbox{in } \ [0,1], \quad \mu=
\frac{q}{q-1}\,,
$$
and $f=0$ at the endpoints, with no free constant. We  get by
integration
 \begin{equation}\label{w3b}
 x-c= z(u,t) \sim (T-t)^{1/(1-q)}g_q(u).
\end{equation}
This formula defines implicitly $u=u(x,t)$ in the space  between
the {\sl converging interfaces}
\begin{equation}\label{w3d}
l(t)=c- C \,(T-t)^{1/(1-q)} , \qquad r(t)=c+  C\, (T-t)^{1/(1-q)}.
\end{equation}
The gradient is given by
\begin{equation}\label{w3c} u_x(x,t)\sim (T-t)^{-1/(1-q)}\frac1{f_q(u)},
\end{equation}
and blows up like $O((T-t)^{1/2\alpha})$ as $t\to T$ and like
$O(d^{-\gamma})$ at the interfaces, $\gamma$ as before.


\section{\textsf{Asymptotic estimates for more general $\Phi$}}
\label{sect-ae2}
\setcounter{equation}{0}

Assumptions on $\Phi$ to behave like a power at infinity guarantee
that the asymptotic behaviour is as predicted in the power case.
Thus, Bertsch and Peletier study in \cite{BP} the asymptotic
behaviour of equation (\ref{eq7}) written in the form
$$
\beta(v)_t=v_{xx}\,.
$$
In our application the space variable $x$ becomes $u$ and  $\beta=\Psi^{-1}$.
They assume that $\beta(0)=0$, $\beta'(s)>0$ for $s>0$ and
$$
-1\le \frac{s\beta''(s)}{\beta'(s)}\le\alpha
$$
for all small $0<s<s_0$ and some $\alpha\in (0,1)$. They also need
the more stringent condition that the limit
$$
\sigma(s)=\lim_{\ve\to 0} \frac{s\beta''(\ve s)}{\beta'(\ve)}
$$
exists for every $s>0$. This limit is necessarily of the form
$\sigma(s)=s^{m}$ for some $m\in [0,1)$. Those conditions cover in
particular powers $\Psi(s)=s^q$ with $q=1/m>1$, as well as the
exponential $\Psi(s)=e^{-1/s}$, and any other function resembling
such examples near $s=0$.

Under these restrictions, they establish
an asymptotic separation-of-variables result that, when translated
to our setting, means that the solution $w(u,t)$ of the Cauchy problem
(\ref{eq7}) behaves as  $t\to \infty$ in the separate variable
form $w\sim y(t)f(u)$, more precisely,
\begin{equation}
\lim_{t\to \infty}\frac{\beta(w)}{y(t)}=f(u),
\end{equation}
where $y(t)$ is the solution $y'=-y^{1/q}$, i.e. $y(t) \sim
c\,t^{-q/(q-1)}$ and $f$ is the profile obtained as in the
previous section. This means that the analysis of the previous
section is justified for {\sl almost power-like} flux functions if
they resemble a power with exponent $q>1$. The cases $q=1$ and
$q<1$ should be justified in a similar manner. We note
that a similar calculation for the Cauchy problem was  done in
\cite{PV}.

We point out that these restrictions on the nonlinearity are
satisfied by the equation proposed in the image processing model,
(\ref{eq3}), hence the results apply to that model.


\section{\textsf{ Numerical experiments}}
\label{sect-ne}

We have computed the solutions of the free boundary problem (${\rm
P}_{\rm II})$) by solving the conjugate problem, which is a
homogeneous Neumann problem, in the cases $\alpha=0$,  $\alpha=1$
and  $\alpha=-1/2$. The conjugate equations are $p$-Laplacian
equations with $p=2, 4$ and $3/2$. The computations have been done
with symmetrical data $x_0(u)=-\cos(\pi u)$ and asymmetrical
$x_0(u)= 2u^4-(1/2)u^6 -1$. The last two examples stabilize in
finite time (since $\alpha<0$).

\

\begin{center}
\includegraphics[width=7cm]{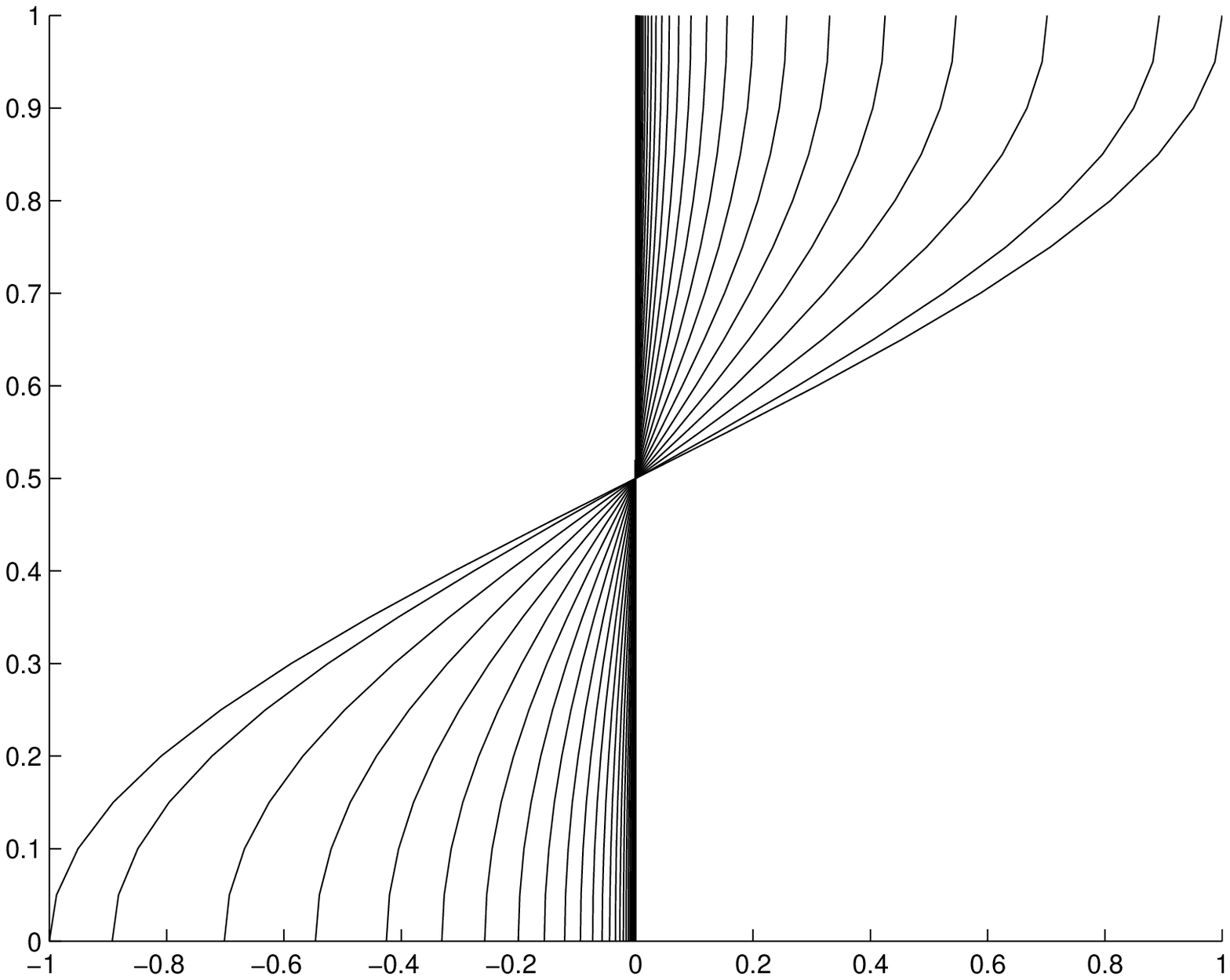}\quad
\includegraphics[width=7cm]{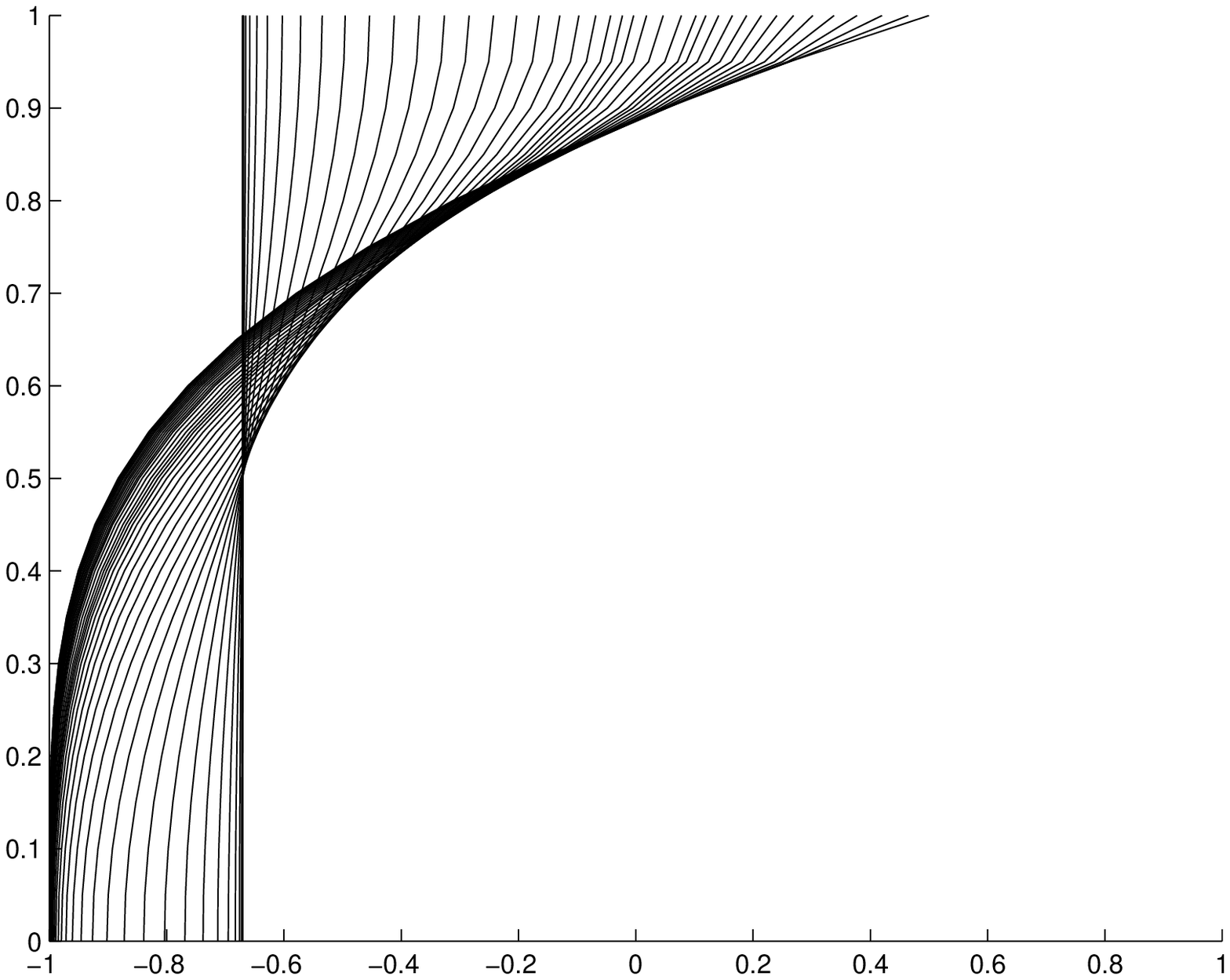}
\end{center}

\centerline{Figures 1 and 2. The case $\alpha=0$ with symmetrical
and asymmetrical data}

\begin{center}
\includegraphics[width=7cm]{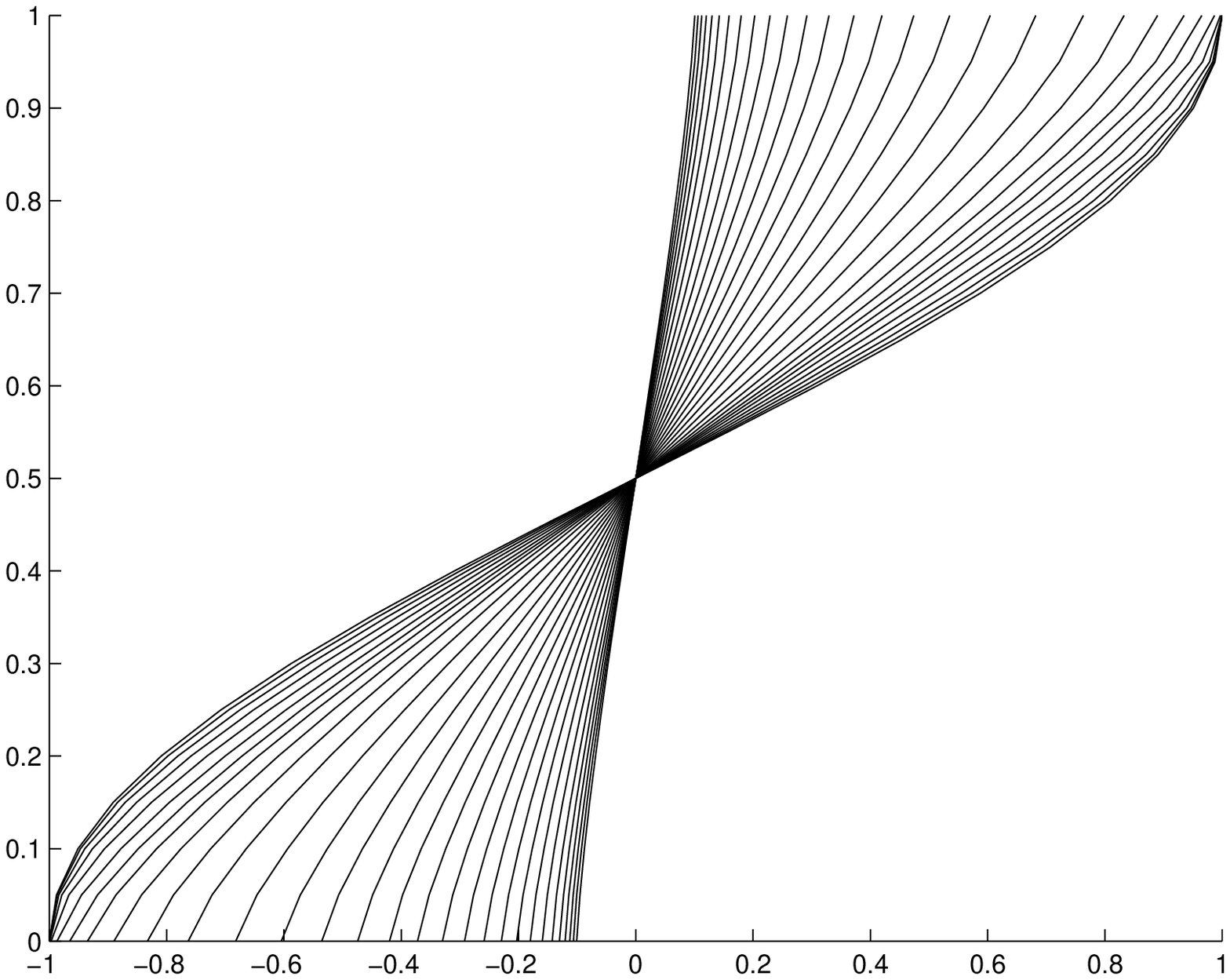}\quad
\includegraphics[width=7cm]{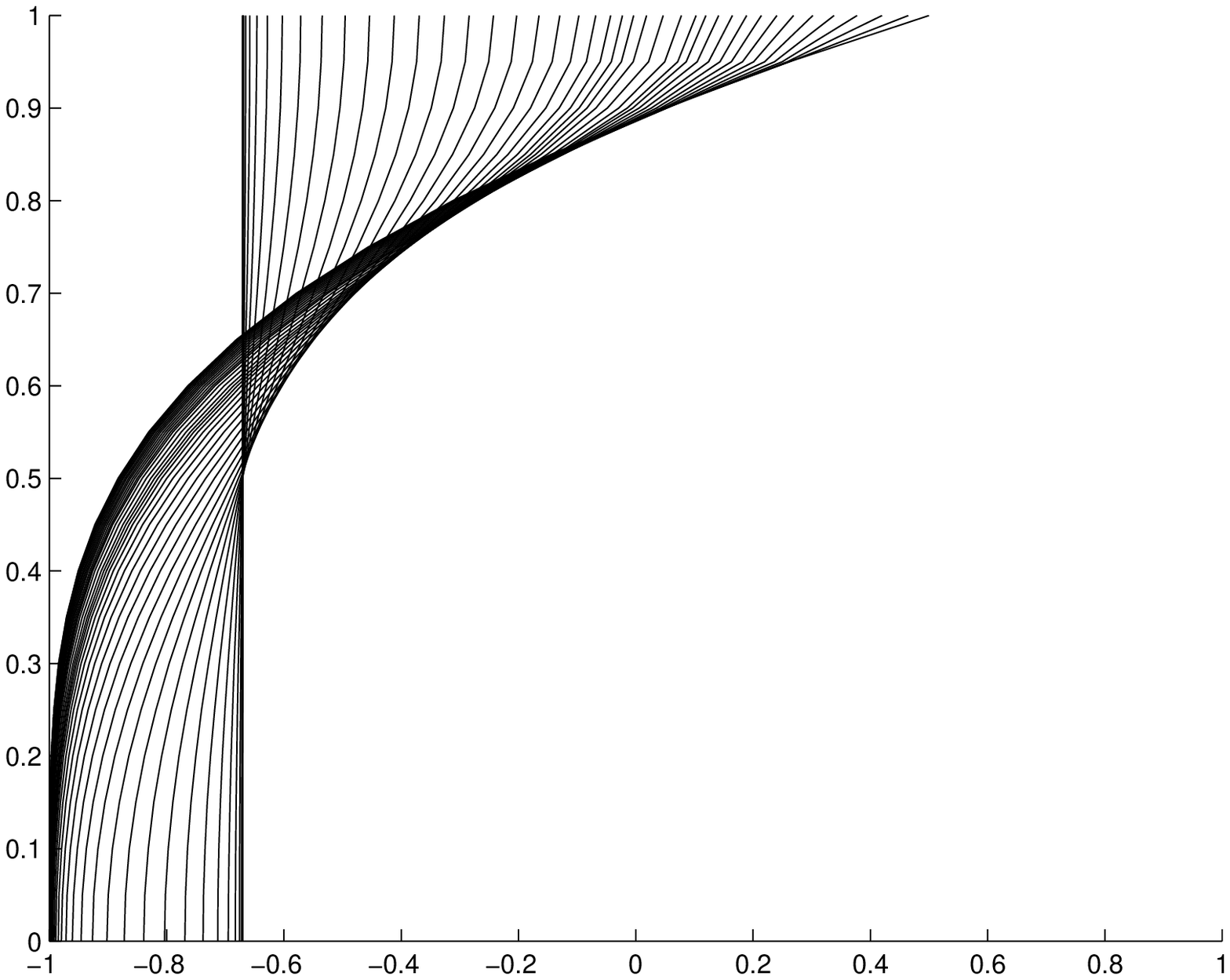}
\end{center}

\centerline{Figures 3 and 4. The case $\alpha=1$ with symmetrical
and asymmetrical data}

\begin{center}
\includegraphics[width=7cm]{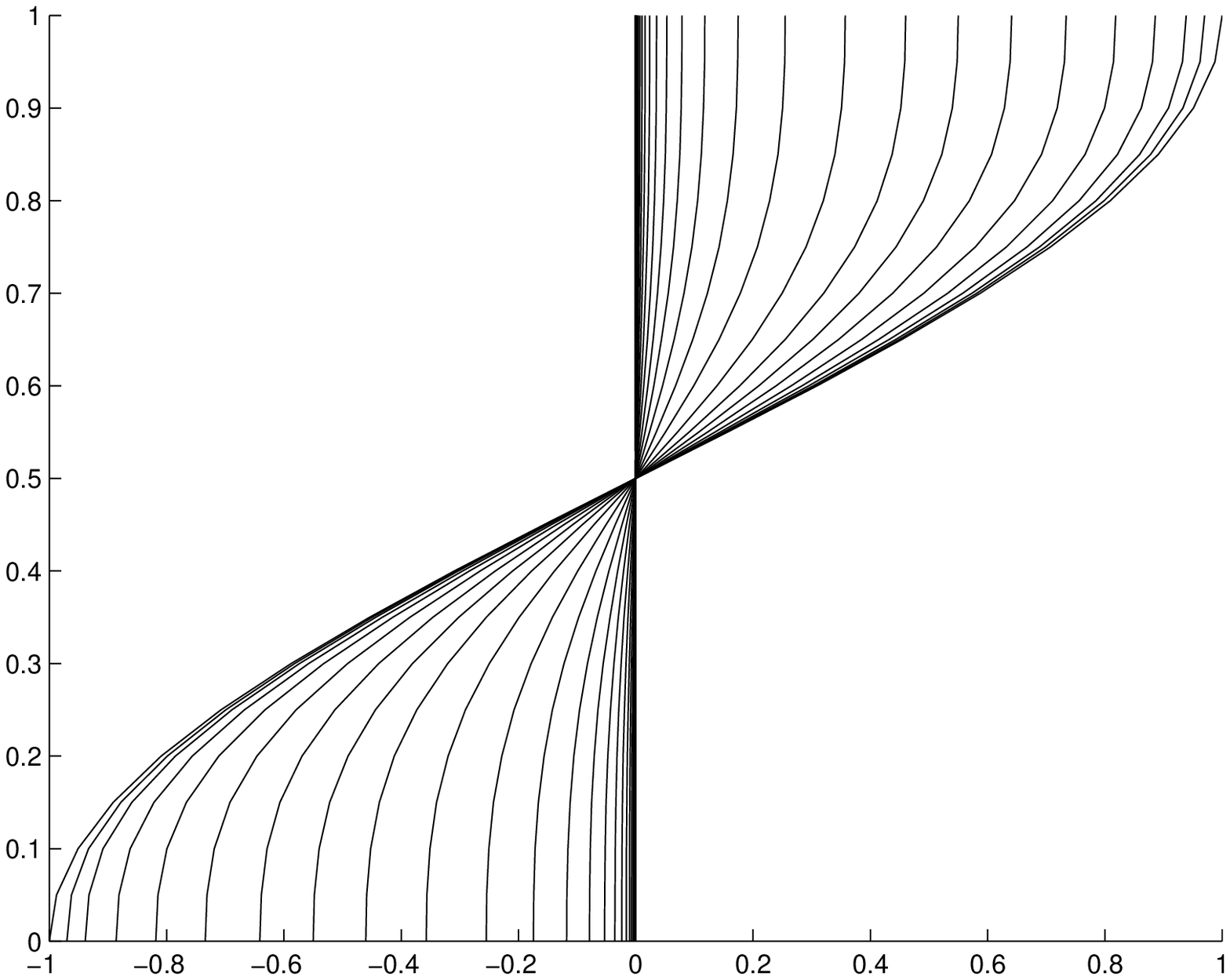}\quad
\includegraphics[width=7cm]{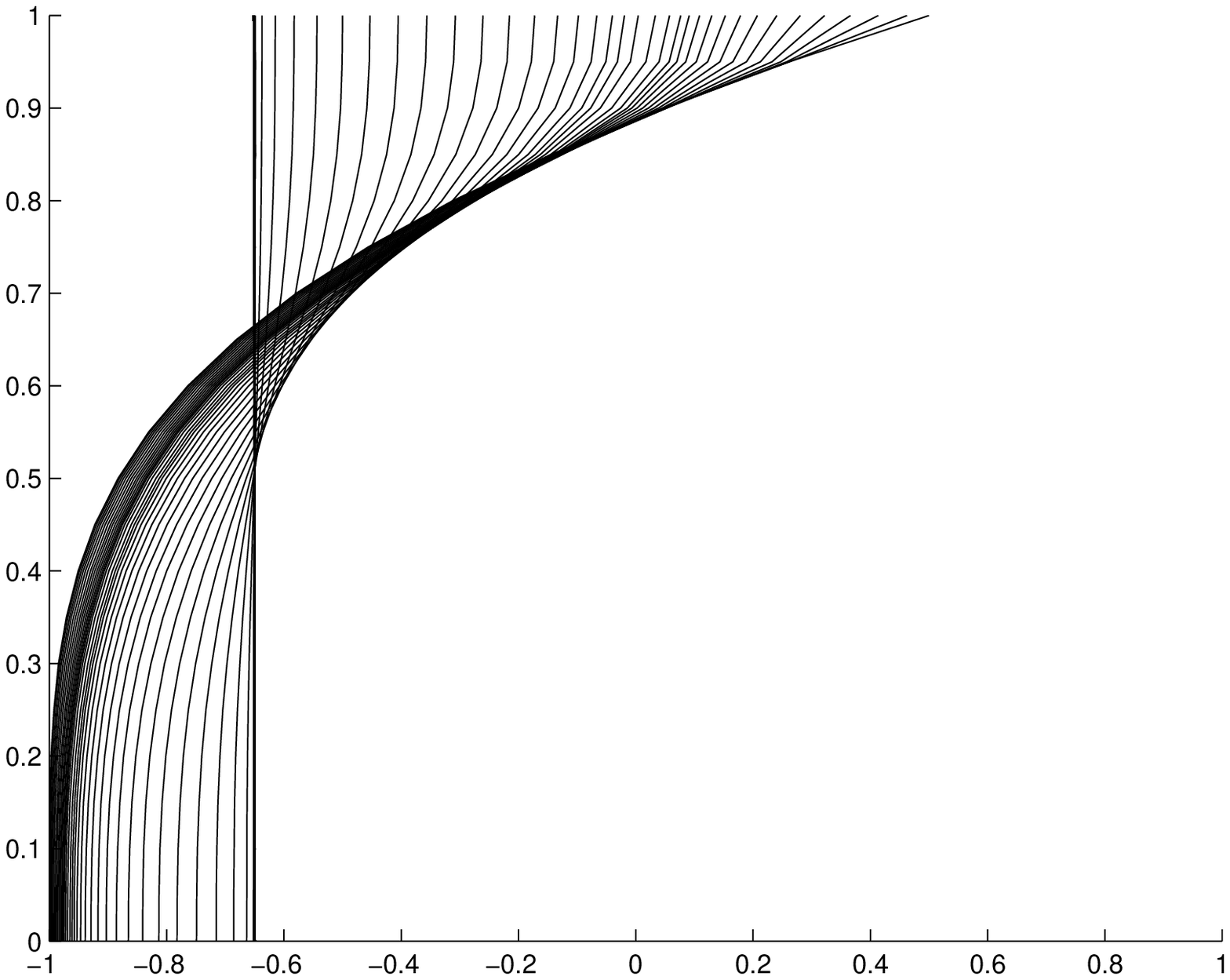}
\end{center}

\centerline{Figures 5 and 6. The case $\alpha=-1/2$ with
symmetrical and asymmetrical data}


\section{\textsf{ Discussion and conclusions }} \label{sect-dc}

We have established the well-posedness of a free boundary problem
that represents a one-dimensional version of the model for image
contour enhancement. We have also established the asymptotic
behaviour and its rates. The results apply to a large class of
equations, which enlarges the scope of the results of \cite{B2,
MS}. We have performed a classification of the solutions and their
properties  according to the properties of the constitutive
function $\Phi$.

 We have used as a framework the theory
nonlinear parabolic equations of diffusive type, or nonlinear
diffusion equations for short. These equations are typically used
in describing processes of mass diffusion or thermal propagation
(with or without additional effects, like convection or reaction).
Here they appear in image processing. The types used here cover
the model cases known as $p$-Laplacian equations and porous medium
equations or their variants.

The analysis is technically performed by means of a series of
remarkable transformations that lead to conjugate problems that
are easier to analyze. These transformations are  related to the
B\"acklund transform.

 There are a number of interesting consequences
of the mathematical analysis that we sum up next.

\noindent $\bullet$ {\sc Generality.} The first one is the observation that
solutions with gradient blowup can be obtained for the nonlinear diffusion equation
with a large class of constitutive functions. Namely, $\Phi$ can be {\sl any
continuous nondecreasing function defined for $s>0$ such that
$\Phi(\infty)$ is bounded}.  We have made for convenience the
assumption of smoothness and strict monotonicity for $0<s<\infty$.
Actually, the constitutive function can be more general, even a
discontinuous maximal monotone graph as in \cite{Be, Brz}. But such a  generality
is not practical here, though the classification and results go through after
some heavy work. A general presentation of that generality will be the object of the
work \cite{Vbk}.

\noindent $\bullet$ {\sc Focusing in finite time.} In the case of
powers $\Phi(s)\sim s^m$ the condition means $m<0$, and not only
$m\le -1$ as considered before. Moreover, the new range $0>m>-1$
leads to a very fast evolution which arrives at a {\sl vertical
front} in finite time. This idea might have an applied interest.

\noindent $\bullet$  {\sc Coexistence of types.} We have also a
so-called class of solutions of Type I  that fall
 into the scope of the standard parabolic theory. They exist under
the condition that the integral $\int (\Phi'(s)/s)\, ds$ converges at
$s=\infty$, which in the case of powers means $m>-1$. This fact leads to an
interesting observation:  there is  a range of co-existence of
solutions of Type I and II for $-1<m<0$, i.e., $-1/2<\alpha<0$.
This may be of interest from the theoretical point of view: {\sl two
different problems share the same equation and same initial data, but the solutions
differ as a consequence of a further choice: of existence or not existence of a free
boundary.}

 \noindent $\bullet$ {\sc The pressure}. There is a way of formally unifying the standard
 theory of free-boundary solutions for  the porous medium equation with
 the  present theory of free-boundaries  in fast diffusion. This
 works by means of the new variable called pressure that for
 an equation of the form $v_t=\Phi(v)_{xx}$ is defined as
\begin{equation}
\label{pr.def}
 p(v)= \int_a^v \frac{\Phi'(s)}{s}\,ds.
 \end{equation}
If the integral is convergent at $v=0$ then the choice $a=0$ is
made. Otherwise, any value in the domain of $\Phi$ is good.
Indeed, the similarity in the behaviour of the pressure of Porous
Medium case, $\Phi(u)=u^m$, $m>1$, and our case (which lies in the
power range $m\le 0$) is seen  when we write the partial
differential equation for $v$ which reads
\begin{equation}
\label{pr.eq} v_t= \sigma(v)\Delta v + |\nabla v|^2,
\end{equation}
with $\sigma(v)=\Phi'(u)$. For $\Phi(u)= cu^m$ we get $v=cmu^{m-1}/(m-1)$
and $\sigma(v)=(m-1)v$. The paper \cite{Vdarcy} explores further the properties of
free-boundary solutions for the pressure of fast diffusion equation. Note that,
contrary to the porous medium case where the support of the pressure solutions
expands in time, in our case it shrinks in time.

\medskip

\noindent $\bullet$ {\sc Grey levels as end values.} Another variation of the main theme is to
consider solutions losing the white or black level. This is
investigated in the paper \cite{RV1}, where many solutions of the
differentiated equation for $v$ are obtained in the range $0\ge
m>-1$ with the same initial data by assigning fluxes at infinity:
\begin{equation}\label{fi}
\lim_{x\to\infty} (v^{m-1}v_x)=-f(t), \qquad \lim_{x\to-\infty}
(v^{m-1}v_x)=g(t),
\end{equation}
for bounded functions $f,g\ge 0$. In terms of $u$ it means in
particular that we can obtain a free-boundary solution with a
decreasing value of highest color level,
$$
u(\infty,t)=F(t)= 1-\int_0^t f(t)\,dt<1.
$$
Analogously, we can impose an increasing value for the lowest
color level,
$$
u(-\infty,t)=G(t)= \int_0^t g(t)\,dt>0.
$$

\noindent $\bullet$ {\sc Non-monotone fronts.} There is a gap in the theory we have developed, namely that we
assumed that the front is monotone in the space direction, in other words,
that $u$ is nondecreasing as a function of $x$. The limitations of fast diffusion equations
to admit non-monotone solutions  are not accidental. In the paper \cite{RV2} it is proved
that the model equation $u_t=u_x^{m-1}u_{xx}$ does not admit non-monotone solutions if
$m\le 0$. In terms of the equation for $v=u_x$ it means that there are no
solutions with both signs. However, there is hope when using functions $\Phi$ like in
(\ref{eq3}) that are degenerate at infinity, i.e., $\Phi'(\infty)=0$,
but regular at all other values, e.g., $\Phi\in C^1(\re)$. Indeed, the problem
with the monotonicity happens because of the singularity at
$u_x=0$ which forbids maxima or minima, while our main interest is in the free-boundary that
is governed by the large values,  $u_x\to\infty$. This question
will be investigated separately.

\noindent $\bullet$ {\sc Non-monotone nonlinearities.} Note that
in our problem setting the equation is forward parabolic, and the
backward movement of the interface is due to the effect of the
singular boundary condition, which happens to be compatible with
the equation for the appropriate class of functions $\Phi$ (e.g.,
for powers $m<0$). In the spirit of the Perona-Malik model,
situations can be considered where $\Phi$ is not monotone, hence
the equation is not parabolic. A regularization is then needed so
that one faces a regularized forward-backward diffusion problem.
See in this respect \cite{BBDU} and \cite{W}. This subject is
again under study.

\


\medskip

 {\sc Acknowledgment}

\noindent The work of the  first author was supported by Applied
Mathematics subprogram of the US Department of Energy under
Contract DE-AC03-76-SF00098. The second author has been partly
supported by  TMR Projects RTN2-2001-00349 and FMRX-CT98-0201, and
Proyecto MCYT BFM2002-04572-C02-02. He is also grateful to  the
Lawrence Berkeley National Laboratory for supporting a research
visit. The authors express their thanks to Prof. A. J. Chorin for
his stimulating interest in our work. Computations were performed
by Ra\'{u}l Ferreira using Matlab.

\setcounter{equation}{0}


\bigskip

\bigskip


{\sc Permanent addresses:}

\textsc{Grigori I. Barenblatt}

Department of Mathematics and Lawrence National Laboratory

Berkeley, California 94720-3840

e-mail: gib@math.uberkeley.edu

\

\textsc{Juan Luis V\'azquez}

Departamento de Matem\'{a}ticas, Universidad Aut\'onoma de Madrid

28046 Madrid, Spain.

e-mail: juanluis.vazquez@uam.es

\bigskip



\begin{thebibliography}{99}





\bibitem[ALM]{ALM}  {\sc L. Alvarez, P.-L. Lions, J.-M. Morel},
Image selective smoothing and edge detection by nonlinear
diffusion. II,
{\it SIAM J. Numer. Anal.} {\bf 29} (1992), 845--866.



\bibitem[ABCM]]{A4} {\sc F. Andreu, C. Ballester, V. Caselles, J.M.
Maz\'on,} Minimizing total variation Flow, {\it Diff. Integral
Eqns.} {\bf 14} (2000), 321--360.



\bibitem[An1]{Ang} {\sc S. B. Angenent}, Local existence and regularity for a
class of degenerate parabolic equations, {\it Math. Ann.}
\textbf{280} (1988), 465--482.



\bibitem[An2]{Ang2} {\sc S. B. Angenent},
 Large time asymptotics for the porous media equation,
{\it  Nonlinear diffusion equations and their equilibrium states,
I} (Berkeley, CA, 1986), 21--34, Math. Sci. Res. Inst. Publ., 12,
Springer, New York, 1988.



\bibitem[Ar]{Ar}  {\sc D. G. Aronson}, {The Porous Medium Equation},
{\it Some problems of Nonlinear Diffusion, Lectures Notes in
Mathematics 1224,} Springer-Verlag, New York, 1986.





\bibitem[AP]{AP}  {\sc D. G. Aronson, L. A. Peletier,}
 Large time behaviour of solutions of the porous medium equation
in bounded domains, {\it J. Diff. Eqns.} {\bf 39} (1981), 378--412.



\bibitem [AV]{AV} {\sc D. G. Aronson,  J.L. Vazquez,}
 Eventual $C^{\infty}$-regularity and concavity of flows in
one-dimensional porous  media,
{\em Archive Rat. Mech. Anal.} {\bf 99} (1987),  329-348



\bibitem[Bbk]{Barbook}  {\sc G. I. Barenblatt},
{\it Scaling, self-similarity, and intermediate asymptotics},
Cambridge Texts in Applied Mathematics, 14. Cambridge University
Press, Cambridge, 1996.



\bibitem[B01]{B2} {\sc  G. I. Barenblatt,}
Self-similar intermediate asymptotics for nonlinear degenerate
parabolic free-boundary problems that occur in image processing,
{\it  Proc. Natl. Acad. Sci. USA} {\bf 98}    (2001), no. 23,
12878--12881.



\bibitem[BBCP] {BBCP} {\sc  G. I. Barenblatt, M. Bertsch, A. E. Chertock, V. M.
Prostokishin}, {Self-similar asymptotics for a degenerate
parabolic filtration-absorption equation}, {\it Proc. Natl. Acad.
Sci. USA} {\bf 97} (2000), no. 18, 9844--9848.





\bibitem[BBDU]{BBDU} {\sc G.I. Barenblatt, M. Bertsch, R. Dal Passo, M.
Ughi,} A degenerate pseudoparabolic regularization of a nonlinear
forward-backward heat equation arising in the theory of heat and
mass exchange in stably stratified turbulent shear flow. {\it SIAM
J. Math. Anal.} {\bf 24} (1993), no. 6, 1414--1439.



\bibitem[BV]{BV} {\sc G. I. Barenblatt, M.I. Vishik},
On finite velocity of propagation in propagation in problems of
non-stationary filtration of a liquid or gas (in Russian), {\it
Prikl. Mat. Mech.} {\bf 20} (1956), 411--417.



\bibitem[Be72]{Be} {\sc P. B\'enilan},  Equations d'\'evolution dans un
espace de Banach quelconque et applications, {\sl Thesis}, Univ.
Orsay, 1972.



\bibitem[Be76] {Be2} {\sc Ph. B\'enilan,}
{Op\'erateurs accr\'etifs et semi-groupes dans les  espaces $L^p$}
($1\le p\le \iy$), in {\it Functional An. and Numerical
Analysis},Japan-France Seminar 1976, ed H. Fujita. Japan Soc.Prom.
Sc. 15--53 (1978).



\bibitem[BC]{BeC} {\sc P. B\'enilan, M. C. Crandall},  The
continuous dependence on \ $\varphi $ \ of Solutions of \ $u_t -
\Delta \varphi (u) = 0$, { \it Indiana Univ. Math. J.}
 { \bf 30} (1981), 161--177.



\bibitem[BG]{BeG} {\sc Ph. B\'enilan, R. Gariepy,}
   {\it Strong solutions in $L^1$ of degenerate parabolic equations},
   J. Diff. Eqns. {\bf 119} (1995), 473--502.



\bibitem[BH]{BH}  {\sc J. G. Berryman, C. J. Holland}, {Stability of the
separable solution for fast diffusion equation,} {\it Arch. Rat. Mech. Anal.}
\textbf{74 }(1980), 379--388.



\bibitem[BPU]{BdPU2} {\sc  M. Bertsch, R. Dal Passo, M. Ughi}, {Nonuniqueness
of solutions of a degenerate parabolic equation,} {\it Annali Mat.
Pura Appl.} \textbf{161} (1992), 57--81.





\bibitem[BP]{BP} {\sc  M. Bertsch, L. A. Peletier},
 The asymptotic profile of solutions of degenerate diffusion
equations, {\it Arch. Rational Mech. Anal.} {\bf 91} (1985), no.
3, 207--229.



\bibitem[BK]{BK} { \sc G. Bluman, S. Kumei},
 On the remarkable nonlinear diffusion equation \break\
$\pa/\pa x [a(u+b)^{-2}(\pa u/\pa x)]-(\pa u/\pa t)=0$,
{\it J. Math. Phys.} { \bf 21} (1980),  1019--1023.



\bibitem[Brz] {Brz} {\sc H. Brezis}, {\it Op\'erateurs maximaux monotones et
semigroupes de contractions dans les espaces de Banach},  (French)
North-Holland Mathematics Studies, No. {\bf 5}. Notas de
   Matem\'{a}tica (50). North-Holland Publishing Co., Amsterdam-London, 1973.





\bibitem[CV]{CV} {\sc L. A. Caffarelli,  J. L. Vazquez}, {Viscosity solutions
for the porous medium equation}, {\it Proc. Symposia in Pure
Mathematics} volume \textbf{65},  Amer. Math. Soc., Providence,
RI, 1999,  13--26. (Proc. of Symposium {\it Differential
equations:} La Pietra 1996 (Florence), in honor of Profs. P. Lax
and L. Nirenberg; M. Giaquinta et al. eds, 1996),



Viscosity solutions for the porous medium
   equation. Differential equations: La Pietra 1996 (Florence), 13--26,
   Proc. Sympos. Pure Math.,
   65, Amer. Math. Soc., Providence, RI, 1999



\bibitem[CMS]{Cas} {\sc V.  Caselles, J. M. Morel, C. Sbert,} An
axiomatic approach to image interpolation.
{\it IEEE Trans. Image Process.} {\bf 7} (1998), no. 3, 376--386.


\bibitem[ChV]{CVp}  {\sc E. Chasseigne, J. L. Vazquez}, The pressure
equation in the Fast Diffusion range, {\it Revista Mat. Iberoamericana},
to appear.



\bibitem[DB]{DiB} {\sc E. DiBenedetto},
{\it Degenerate Parabolic Equations}, Universitext, Springer-
Verlag (1994).



\bibitem[ERV]{ERV} {\sc J. R. Esteban, A. Rodr\'{\i}guez, J. L. Vazquez}, {A
nonlinear heat equation with singular diffusivity}, {\it Comm.
Partial Diff. Eqs.} \textbf{13} (1988), 985-1039.



\bibitem[FV]{FV} {\sc R. Ferreira, J. L. Vazquez},
{Study of self-similarity for the fast diffusion equation},
{\em Advances Diff. Eqns.}, to appear.





\bibitem[FK]{FrKa80}
{\sc  A. Friedman, S. Kamin.}
\newblock The asymptotic behavior of gas in an N-dimensional
porous medium,
\newblock {\em Trans. Amer. Math. Soc.} {\bf 262} (1980), 551--563.



\bibitem[GP]{GP} {\sc B. H. Gilding, L. A. Peletier,} On a class of similarity
solutions of the porous media equation, {\it J. Math. Anal. Appl.} {\bf 55} (1976),
351--364.



\bibitem[Ka] {Kal}  {\sc A. S. Kalashnikov,}
Some problems of the qualitative theory of non-linear
 degenerate second-order parabolic equations,
{ \it Russian Math. Surveys,} {\bf 42} (1987), 169--222.



\bibitem[Kch]{Kc} {\sc S. Kichenassamy}, The Perona-Malik paradox,
{\it SIAM J. Appl. Math.} {\bf 57} (1997), 1328--1324.





\bibitem[Ki]{Ki} {\sc J. R. King},
{Self-similar behaviour for the equation of fast nonlinear
diffusion}, {\it Phil. Trans. Roy. Soc. London A} \textbf{343}
(1993), 337--375.





\bibitem[KB]{KB} {\sc H. Krim, Y. Bao},
{Nonlinear diffusion: a probablilistic approach},
{\em ICP99,} Kobe, Japan.



\bibitem[MS]{MS} {\sc R. Malladi,  J.A. Sethian,}
Image Processing: Flows under Min/Max curvature and Mean
Curvature,
{\it Graphical Models and Image Processing,} {\bf 58}, 2 (1996),
127--141.



\bibitem[OKC]{OKC} {\sc O. Ole\u{\i}nik, S. A. Kalasknikov, Y.-L.
Czhou,}
The Cauchy problem and boundary-value problems for equations of
the type of unsteady filtration,
{\em Izv. Akad. Nauk SSSR,} Ser. Mat. {\bf 22} (1958),  667--704.



\bibitem[PV]{PV} {\sc A. de Pablo, J. L. Vazquez},
 Regularity
of solutions and interfaces of a generalized porous medium
equation, {\it Ann. Mat. Pura Applic. (IV),} {\bf 158} (1991),
51--74.



\bibitem[PM]{PM} {\sc P. Perona, J. Malik},
Scale space and edge detection using anisotropic diffusion,
{\it IEEE Transactions of Pattern Analysis and Machine
Intelligence}, {\bf 12} (1990), 629--639.



\bibitem[RV93] {RV93} {\sc A. Rodriguez, J. L. Vazquez,}
{\it  A well-posed problem in singular Fickian diffusion},
Archive Rat. Mech. Anal. {\bf 110}, 2 (1990),
141--163.



\bibitem[RV95] {RV1} {\sc A. Rodriguez, J. L. Vazquez,}
{\it Non-uniqueness of solutions of nonlinear heat equations of
fast diffusion type}, Ann. Inst. Henri Poincar\'e, Analyse Non
Lin., 12, 2 (1995), 173--200.



\bibitem[RV02]{RV2} {\sc A. Rodr\'{\i}guez, J. L. Vazquez}, {Obstructions to
existence in fast-diffusion equations}, {\it  J. Diff. Eqns.} {\bf 184}
(2002), 348--385.





\bibitem[SGKM]{S4} {\sc A. A. Samarskii, V. A. Galaktionov,
S. P. Kurdyumov,  A. P. Mikhailov}, \lq \lq Blow-up in problems
for quasilinear parabolic equations''. Nauka, Moscow, 1987 (in
Russian). English transl.: Walter de Gruyter, Berlin, 1995.



\bibitem[SKMS]{SKM} {\sc N. Sochen, R. Kimmel, R. Malladi,}
From high energy physics to low level vision, {\it Preprint
LBNL-39243, UC-405}, E. O. Lawrence Berkeley National Laboratory
(1996), 38 pags.





\bibitem[SKM]{SKM2} {\sc N. Sochen, R. Kimmel, R. Malladi},
A general framework for low level vision, {\it IEEE Trans. Image
Process.} {\bf 7} (1998), 310--318.



\bibitem[V92]{Vb} {\sc J. L. Vazquez}, ``An introduction to the mathematical
theory of the porous medium equation'', in {\it  Shape
optimization and free boundaries} (Montr\'eal, PQ, 1990),
347--389, NATO Adv. Sci. Inst. Ser. C Math. Phys. Sci., 380,
Kluwer Acad. Publ., Dordrecht, 1992.



\bibitem[V03a]{Vncp} {\sc J. L. Vazquez},
 Asymptotic behaviour for the PME in the whole space, \newline
 {\em  Journal of Evolution Equations}, to appear. UAM Preprint
 {\tt http://www.uam.es/ juanluis.vazquez}.



\bibitem[V03b]{V3b} {\sc J. L. Vazquez},
``Asymptotic behaviour for the PME in a bounded domain. The
Dirichlet problem". {\em Monatshefte f. Mathematik}, to appear.
Course Notes, UAM,  1997, {\tt http://www.uam.es/
juanluis.vazquez.}



\bibitem[V03c]{Vdarcy} {\sc J.L. Vazquez}, Darcy's Law and the theory
of shrinking solutions of fast diffusion equations,  {\it TICAM
Report 01-18}, Univ. of Texas at Austin, 2001.



\bibitem[V03d]{Vposit} {\sc J.L. Vazquez},
{Positivity, Propagation Properties and Needles in Nonlinear
Singular Diffusion}, {\em Preprint Univ. Aut\'onoma Madrid.}


\bibitem[Vbk]{Vbk} {\sc J.L. Vazquez}, {\it One-Dimensional Nonlinear
Diffusion},  Course Notes, UAM 1999. Book in preparation.



\bibitem[Wck]{Weick} {\sc J. Weickert}, Applications of nonlinear diffusion in
image processing and computer vision. {\sl Acta Math. Univ. Comenian.}
(N.S.) 70 (2000), no. 1, 33--50.



\bibitem[WB]{WB} {\sc J. Weickert, B. Benhamouda},
A semidiscrete nonlinear scale-space theory and its relation to the Perona-Malik
paradox. F. Solina, W. G. Kropatsch, R. Klette, R. Bajcsy (Eds.),
{\em Theoretical foundations of computer vision (TFCV '96)}, Dagstuhl,
March 18-22, 1996), Springer, Wien, pp. 1-10, 1997.



\bibitem[Wit]{W} {\sc T. P. Witelski,}
The structure of internal layers for unstable nonlinear diffusion
equations. {\it Stud. Appl. Math.} {\bf 97} (1996), no. 3,
277--300.



\end{thebibliography}
\end{document}